\documentclass[11pt]{amsart}
\usepackage{amsmath,amsthm,amssymb,amsrefs}
\usepackage{tikz-cd}
\usepackage{enumerate}
\usetikzlibrary{fit,shapes.geometric}
\newtheorem{thm}{Theorem}[section]
\newtheorem{cor}[thm]{Corollary}

\newtheorem{prp}[thm]{Proposition}

\newtheorem{con}[thm]{Conjecture}
\theoremstyle{definition}
\newtheorem{df}[thm]{Definition}

\newtheorem{prb}[thm]{Problem}

\theoremstyle{remark}

\numberwithin{equation}{section}
\newcommand{\meet}{\wedge}

\begin{document}

\title[The Complexity of Homomorphism Factorization]
{The Complexity of Homomorphism Factorization}
\author[Kevin M. Berg]
{Kevin M. Berg}
\address[Kevin M. Berg]{
Department of Mathematics\\
University of Colorado\\
Boulder, CO 80309-0395, U.S.A.}
\email{Kevin.Berg@Colorado.EDU}
\thanks{This material is based upon work supported by the National Science Foundation grant no.\ DMS 1500254.}
\subjclass{Primary 68Q25, Secondary 03C05, 08A05}
\keywords{homomorphism factorization, NP-complete, computational complexity, finite structures}

\begin{abstract} 
	We investigate the computational complexity of the problem of deciding if an algebra homomorphism can be factored through an intermediate algebra. Specifically, we fix an algebraic language, $\mathcal L$, and take as input an algebra homomorphism $f\colon X\to Z$ between two finite $\mathcal L$-algebras $X$ and $Z$, along with an intermediate finite $\mathcal L$-algebra $Y$. The decision problem asks whether there are homomorphisms $g\colon X\to Y$ and $h\colon Y\to Z$ such that $f=hg$. We show that these Homomorphism Factorization Problems are NP-complete. We also develop a technique for producing compatible restrictions on homomorphisms, and show that Homomorphism Factorization Problems have polynomial time instances for finite Boolean algebras, finite vector spaces, finite $G$-sets, and finite abelian groups.
\end{abstract}

\maketitle

\section{Introduction}\label{sec:Intro}

\begin{figure}[h]
			\begin{tikzcd}[sep=large]
				X \arrow[dr, dotted, "{\exists g?}"'] \arrow[rr, "f"] && Z\\
				& Y \arrow[ur, dotted, "{\exists h?}"']
			\end{tikzcd}
		\caption{The general form of the commutative diagram for Homomorphism Factorization Problems.}
		\label{fig:diagram}
	\end{figure}
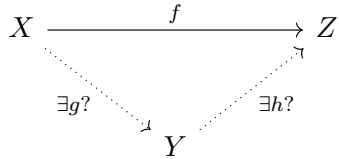

	In this paper we investigate the complexity of the problem of deciding if an algebra homomorphism can be factored through an intermediate algebra. Specifically, we fix an algebraic language $\mathcal L$. The input to our problem is a homomorphism $f\colon X\to Z$ between $\mathcal L$-algebras $X$ and $Z$, along with an intermediate $\mathcal L$-algebra $Y$. The problem is to decide whether there are homomorphisms $g\colon X\to Y$ and $h\colon Y\to Z$ such that $f=hg$, as shown in Fig. \ref{fig:diagram}. We refer to this as
	
	\begin{prb}[The Homomorphism Factorization Problem]\label{HFP}
		Given a homomorphism $f\colon X\to Z$ between two finite $\mathcal L$-algebras $X$ and $Z$, and given an intermediate finite $\mathcal L$-algebra $Y$, decide whether there are homomorphisms $g\colon X\to Y$ and $h\colon Y\to Z$ such that $f=hg$.
	\end{prb}

	There are several interesting special cases of the main problem worth identifying:

\begin{enumerate}[I.]
	\item Consider restricting to the case where $|Z|=1$. Then the homomorphisms $f$ and $h$ from Problem \ref{HFP} must be constant, so the Homomorphism Factorization Problem reduces to the problem of deciding whether, given $\mathcal L$-algebras $X$ and $Y$, there is a homomorphism $g\colon X\to Y$, as shown in Fig. \ref{fig:hom}. We refer to this special case as the \emph{Homomorphism Problem}.
				
	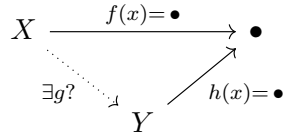
\begin{figure}[h]
			\begin{tikzcd}
				X \arrow[dr, dotted, "\exists g?"'] \arrow[rr, "{f(x) = \, \bullet}"] && \bullet \\
				& Y \arrow[ur, "{h(x) = \, \bullet}"']
			\end{tikzcd}
		\caption{The commutative diagram for the Homomorphism Problem.}
		\label{fig:hom}
	\end{figure}

	\item Consider restricting to the case where the input is $\mathcal L$-algebras $X$, $Y$, $Z$, and homomorphisms $f\colon X \to Z$, and $h\colon Y \to Z$, so the Homomorphism Factorization Problem reduces to the problem of deciding whether there is a homomorphism $g\colon X \to Y$, as shown in Fig. \ref{fig:erfp}. This special case will be called the \emph{Exists Right-Factor Problem}.\\
	
	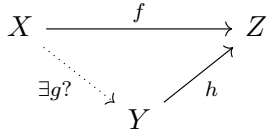
\begin{figure}[h]
			\begin{tikzcd}
				X \arrow[dr, dotted, "\exists g?"'] \arrow[rr, "f"] && Z \\
				& Y \arrow[ur, "h"']
			\end{tikzcd}
		\caption{The commutative diagram for the Exists Right-Factor Problem.}
		\label{fig:erfp}
	\end{figure}
	
	\item Consider restricting to the case where the input is $\mathcal L$-algebras $X$, $Y$, $Z$, and homomorphisms $f\colon X \to Z$, and $g\colon X \to Z$, so the Homomorphism Factorization Problem reduces to the problem of deciding whether there is a homomorphism $h\colon Y \to Z$, as shown in Fig. \ref{fig:elfp}. This special case will be called the \emph{Exists Left-Factor Problem}.\\
	
	\begin{figure}[h]
			\begin{tikzcd}
				X \arrow[dr, "g"'] \arrow[rr, "f"] && Z \\
				& Y \arrow[ur, dotted, "\exists h?"']
			\end{tikzcd}
		\caption{The commutative diagram for the Exists Left-Factor Problem.}
		\label{fig:elfp}
	\end{figure}
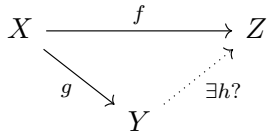

	\item Consider restricting to the case where $f\colon X\to Z$ is the identity function from $X$ to $Z=X$, as shown in Fig. \ref{fig:ret}. In this case the Homomorphism Factorization Problem reduces to the problem of deciding if, given $X$ and $Y$, the algebra $X$ is a retract of $Y$. This special case will be called the \emph{Retraction Problem}.\\
	
	\begin{figure}[h]
			\begin{tikzcd}
				X \arrow[dr, dotted, "\exists g?"'] \arrow[rr, "id"] && X \\
				& Y \arrow[ur, dotted, "\exists h?"']
			\end{tikzcd}
		\caption{The commutative diagram for the Retraction Problem.}
		\label{fig:ret}
	\end{figure}
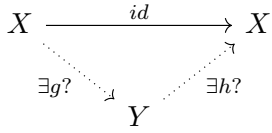	

	\item Consider restricting the Retraction Problem in the case where $|X|=|Y|$. This special case is the \emph{Isomorphism Problem} for $\mathcal L$-algebras.
\end{enumerate}

	The problem of deciding the complexity of the Homomorphism Factorization Problem was raised at \cite{JosephVanName}. There the focus is on algebras in the language of a single fundamental operation, which is binary. The author of the problem expresses interest in the special case where the algebras are semigroups. We shall also focus on these special cases, and we shall consider only finite algebras.

	The questions considered here have been studied for relational structures, so it is important to indicate the differences. Any semigroup can be viewed are a relational structure with a single ternary relation $\{(x,y,z)\in X^3\;|\;z=xy\}$. A function $f\colon X\to Z$ between semigroups is an algebra homomorphism when $X$ and $Z$ are considered as algebras if and only if it is a relational homomorphism when $X$ and $Z$ are considered as relational structures. Therefore, the problem of deciding if a semigroup algebra homomorphism can be factored is the same as the problem of deciding if a semigroup relational homomorphism can be factored.

	But the problem of deciding if a semigroup homomorphism can be factored is not the same as the problem of deciding if a homomorphism of relational structures (with one ternary relation) can be factored. The latter problem involves relational structures that are not codings of semigroups. In fact, it is not hard to see, and it is well known, that the homomorphism problem for ternary relational structures is NP-complete. But the homomorphism problem for semigroups always has an affirmative answer. That is, given finite semigroups $X$ and $Y$, there is always a semigroup homomorphism $g\colon X\to Y$, namely we can take $g$ to be a constant homomorphism mapping $X$ to an idempotent of $Y$.

	It is also worth noting that some cases of the Homomorphism Factorization Problem can be easy (e.g., we just noted that the Homomorphism Problem for semigroups always has an affirmative answer), but also can be hard (e.g., the Group Isomorphism Problem is a special case, and there is no known easy algorithm to decide the Group Isomorphism Problem). In order to characterize these cases, we introduce the notion of a ``rich'' language. An algebraic language is \emph{rich} if it has at least one operation of arity at least two, or at least two unary operations. We will prove the following:
	
	\begin{thm}\label{thm:General}The Homomorphism Factorization Problem is NP-complete for rich languages.
	\end{thm}

\section{Rich Languages with at Least Two Unary Operations}\label{sec:ThreeUnary}
	
	Recall that an algebraic language is \emph{rich} if it has at least one operation of arity at least two, or at least two unary operations. In this section, we begin the proof of Theorem \ref{thm:General} by considering the Exists Right-Factor Problem for algebras with two or more unary operations. Specifically, we will show that in this case, any algorithmic solution to the Homomorphism Problem for such algebras would necessarily give a solution to the Homomorphism Problem for graphs, and vice versa. As Graph Homomorphism is well-known to be NP-complete \cite{GareyJohn}, this will provide the desired result.
	
	We may consider any loop-free, finite, connected, directed graph, $G = (V_G, E_G)$, with at least two vertices. We encode this graph into an algebra $G^{\dagger}$, using the following rules to construct the universe: for each vertex $v$ in $V_G$, there are two corresponding elements $v_1$ and $v_2$ in $G^{\dagger}$, and for each edge $(u,v)$ in $E_G$, there are two elements, $a_{(u,v)}$ and $b_{(u,v)}$, in $G^{\dagger}$. We assign to $G^{\dagger}$ two unary operations: $f(\cdot)$ and $g(\cdot)$, given by Table \ref{tab:unary}.
	
	\begin{table}[h]
 		\caption{The operations of $G^{\dagger}$.}
		\label{tab:unary}
		\begin{tabular}{l|ll}
			${}$   & $f(\cdot)$ & $g(\cdot)$ \\ \hline
 			$u_1$ & $u_1$ & $u_2$ \\
 			$u_2$ & $u_1$ & $u_2$ \\
 			$a_{(u,v)}$ & $u_1$ & $b_{(u,v)}$ \\
 			$b_{(u,v)}$ & $v_2$ & $a_{(u,v)}$
 		\end{tabular}
	\end{table}
	
	This encoding allows us to prove the following:
	
	\begin{thm} Let $G$ and $H$ be loop-free, finite, connected, directed graphs with at least two vertices. There exists a homomorphism $\phi \colon G \to H$ if and only if there exists a homomorphism $\psi \colon G^{\dagger} \to H^{\dagger}$.
	
		\begin{proof} Suppose first that there exists a homomorphism $\phi \colon G \to H$. We construct a function $\psi\colon G^{\dagger} \to H^{\dagger}$ based on $\phi$ -- specifically, if for any $v$ in $V_G$ we have $\phi(v)$ in $V_H$, then $\psi(v_1) = \phi(v)_1$, $\psi(v_2) = \phi(v)_2$ and for all $(u,v)$ in $E_G$, $\psi(a_{(u,v)}) = a_{(\phi(u),\phi(v))}$ and $\psi(b_{(u,v)}) = b_{(\phi(u),\phi(v))}$.
		
		We claim $\psi$ is well-defined by the well-definition of $\phi$ and the construction of $H^{\dagger}$. Since $u$ maps to a single element of $H$ under $\phi$, that single element, $\phi(u)$, is in turn associated with exactly two distinct elements $\phi(u)_1$ and $\phi(u)_2$ in $H^{\dagger}$. Similarly, for any $(u,v)$ in $E_G$, we must have that $(\phi(u),\phi(v))$ is in $E_H$. In turn, it must be the case that $a_{(\phi(u),\phi(v))}$ and $b_{(\phi(u),\phi(v))}$ are two single elements in $H^{\dagger}$.
		
		Next, we claim $\psi$ is a homomorphism. Suppose $i$ and $j$ are distinct elements of $\{1,2\}$, and $v_i$ is an element of $G^{\dagger}$ coming from any $v$ in $V_G$. We have that $\psi(f(v_i)) = \psi(v_1) = \phi(v)_1 = f(\phi(v)_i) = f(\psi(v_i))$, and $\psi(g(v_i)) = \psi(v_2) = \phi(v)_2 = g(\phi(v)_i) = g(\psi(v_i))$. Similarly, if $a_{(u,v)}$ and $b_{(u,v)}$ are elements of $G^{\dagger}$, then we have that $\psi(f(a_{(u,v)})) = \psi(u_1) = \phi(u)_1 = f(a_{(\phi(u),\phi(v))}) = f(\psi(a_{(u,v)}))$, $\psi(g(a_{(u,v)})) = \psi(b_{(u,v)}) = b_{(\phi(u),\phi(v))} = g(a_{(\phi(u),\phi(v))}) = g(\psi(a_{(u,v)}))$, $\psi(f(b_{(u,v)})) = \psi(v_2) = \phi(v)_2 = f(b_{(\phi(u),\phi(v))}) = f(\psi(b_{(u,v)}))$, and $\psi(g(b_{(u,v)})) = \psi(a_{(u,v)}) = a_{(\phi(u),\phi(v))} = g(b_{(\phi(u),\phi(v))}) = g(\psi(b_{(u,v)}))$.
		
		Suppose now that there instead exists a homomorphism $\psi\colon G^{\dagger} \to H^{\dagger}$.
		
		We claim $\psi(u_2) = u'_2$ for some $u' \in V_H$. By construction, $u_2 = gf(u_2)$, and since $\psi$ is a homomorphism, $\psi(u_2) = \psi(gf(u_2)) = gf(\psi(u_2))$. The only $x \in H^{\dagger}$ satisfying $x = gf(x)$ are $x = u'_2$ for some $u' \in V_H$. The claim holds. Furthermore, since $\psi(u_2) = u'_2$ for some $u' \in V_H$, it follows that $\psi(u_1) = u'_1$ for the same $u' \in V_H$ -- since $f(u_2) = u_1$, $\psi(u_1) = \psi(f(u_2)) = f(\psi(u_2)) = f(u'_2) = u'_1$.
		
		We now claim that $\psi(b_{(u,v)}) = b_{(u',v')}$ for some $u', v' \in V_H$. Since $f(b_{(u,v)}) = v_2 = g(f(b_{(u,v)}))$, $f(\psi(b_{(u,v)}) = \psi(f(b_{(u,v)})) = \psi(g(f(b_{(u,v)}))) = g(f(\psi(b_{(u,v)})))$. Then $\psi(b_{(u,v)})$ satisfies the equational condition that $f(x) = gf(x)$, and the only $x$ with $f(x) = gf(x)$ are precisely the elements of the form $b_{(u',v')}$ for some $u', v' \in V_H$. The claim holds. Furthermore, since $\psi(b_{(u,v)}) = b_{(u',v')}$ for some $u', v' \in V_H$, it follows that $\psi(a_{(u,v)}) = a_{(u',v')}$ for the same $u', v' \in V_H$ -- since $g(b_{(u,v)}) = a_{(u,v)}$, $\psi(a_{(u,v)}) = \psi(g(b_{(u,v)})) = g(\psi(b_{(u,v)})) = g(b_{(u',v')}) = a_{(u',v')}$.
		
		This provides the basis for our construction of $\phi$ -- if $v$ is in $V_G$, we let $\phi(v) = w$ in $V_H$ where $\psi(v_1) = w_1$. By the preceding arguments, $\phi$ is a well-defined homomorphism.
		\end{proof}
		
	\end{thm}
	
	\begin{cor}\label{cor:unary} The Homomorphism Problem for algebras in a language with at least two unary operations is NP-complete.\qed\end{cor}
	
	In fact, this same construction allows us to completely classify the Homomorphism Factorization Problem for such languages. The observation that the Homomorphism Problem is a special case of the Exists Right-Factor Problem gives the following:
	
	\begin{cor} The Exists-Right Factor Problem for algebras in a language with at least two unary operations is NP-complete.\qed\end{cor}
	
	We can also classify the Exists-Left Factor Problem using a specific input corresponding to a homomorphism with two fixed points.
	
	\begin{cor} The Exists-Left Factor Problem for algebras in a language with at least two unary operations is NP-complete.
		
		\begin{proof} Let $G=(\{v\}, \emptyset)$ be the graph consisting of a single vertex, $v$, with no edges. We encode $G$ into an algebra, $X$, using the standard rules for constructing $G^{\dagger}$ -- specifically, the universe of $X$ consists of two elements, $v_1$ and $v_2$, and we assign to $X$ two unary operations, $f(\cdot)$ and $g(\cdot)$, given by Table \ref{tab:X}.
		
			\begin{table}[h]
 				\caption{The operations of $X$.}
				\label{tab:X}
				\begin{tabular}{l|ll}
					${}$   & $f(\cdot)$ & $g(\cdot)$ \\ \hline
 					$v_1$ & $v_1$ & $v_2$ \\
 					$v_2$ & $v_1$ & $v_2$
 				\end{tabular}
			\end{table}
			
			We then set algebras $Y$ and $Z$ to be the encodings of two loop-free, finite, connected, directed graphs with at least two vertices, $H$ and $J$, each together with a new distinguished vertex $v'$ not connected to any other vertices. In particular, $v'$ is associated with neither $a$ nor $b$ elements in either algebra. Suppose $f \colon X \to Z$ is given by $f(v) = v'$ and $g\colon X \to Y$ is given by $g(v) = v'$. Then a homomorphism $h\colon Y \to Z$ with $f = hg$ exists if and only if there exists a homomorphism between $H^{\dagger}$ and $J^{\dagger}$. By Corollary \ref{cor:unary}, the determination of the existence of such an $h$ is NP-complete.
		\end{proof}
	
	\end{cor}
	
	By the NP-completeness of the Subgraph Isomorphism Problem \cite{GareyJohn}, this construction also yields the following pair of results:
	
	\begin{thm} Let $G$ and $H$ be loop-free, finite, connected, directed graphs with at least two vertices. Then $G$ is isomorphic to a subgraph of $H$ if and only if $G^{\dagger}$ is a retraction of $H^{\dagger}$.
	
		\begin{proof} Suppose first that $G$ is isomorphic to a subgraph of $H$. Then, in particular, there exist graph homomorphisms $\phi_1 \colon G \to H$ and $\phi_2 \colon H \to G$ with $\phi_1\phi_2 = \text{id}_H$, the identity map on $H$. By Theorem \ref{thm:unary}, there must exist homomorphisms $\psi_1\colon G^{\dagger} \to H^{\dagger}$ and $\psi_2\colon H^{\dagger} \to G^{\dagger}$. In particular, for any $v$ in $V_G$, we must have $\psi_1(v_1) = \phi_1(v)_1$ and $\psi_1(v_2) = \phi_1(v)_2$. Similarly, for any $(u,v)$ in $E_G$, we must have $\psi_1(a_{(u,v)}) = a_{(\phi_1(u),\phi_1(v))}$ and $\psi_1(b_{(u,v)}) = b_{(\phi_1(u),\phi_1(v))}$. Correspondingly, for any $v$ in $V_H$ or $(u,v)$ in $E_H$, we have $\psi_2(v_1) = \phi_2(v)_1$, $\psi_2(v_2) = \phi_2(v)_2$, $\psi_2(a_{(u,v)}) = a_{(\phi_2(u),\phi_2(v))}$, and $\psi_2(b_{(u,v)}) = b_{(\phi_2(u),\phi_2(v))}$.
		
		We claim that $\text{id}_{H^{\dagger}}$, the identity map on $H^{\dagger}$, must be equal to $\psi_1\psi_2$. Let $v$ be any vertex in $V_H$. Then $\text{id}_{H^{\dagger}}(v_1) = v_1 = \text{id}_H(v)_1 = \phi_1(\phi_2(v))_1 = \psi_1(\phi_2(v)_1)=\psi_1(\psi_2(v_1)) = \psi_1\psi_2(v_1)$ and $\text{id}_{H^{\dagger}}(v_2) = v_2 = \text{id}_H(v)_2 = \phi_1(\phi_2(v))_2 = \psi_1(\phi_2(v)_2)=\psi_1(\psi_2(v_2)) = \psi_1\psi_2(v_2)$. Similarly, let $(u,v)$ be any element of $E_H$. Then $\text{id}_{H^{\dagger}}(a_{(u,v)}) = a_{(u,v)} = a_{(\text{id}_H(u),\text{id}_H(v))} = a_{(\phi_1\phi_2(u),\phi_1\phi_2(v))} = \psi_1(a_{(phi_2(u),\phi_2(v))})=\psi_1(\psi_2(a_{(u,v)})) = \psi_1\psi_2(a_{(u,v)})$ and $\text{id}_{H^{\dagger}}(b_{(u,v)}) = b_{(u,v)} = b_{(\text{id}_H(u),\text{id}_H(v))} = b_{(\phi_1\phi_2(u),\phi_1\phi_2(v))} = \psi_1(b_{(phi_2(u),\phi_2(v))})=\psi_1(\psi_2(b_{(u,v)})) = \psi_1\psi_2(b_{(u,v)})$. The claim holds.
		
		Suppose now that $G^{\dagger}$ is a retraction of $H^{\dagger}$. Then, in particular, there exist homomorphisms $\psi_1 \colon G^{\dagger} \to H^{\dagger}$ and $\psi_2 \colon H^{\dagger} \to G^{\dagger}$ with $\psi_1\psi_2 = \text{id}_{H^{\dagger}}$, the identity map on $H^{\dagger}$. By Theorem \ref{thm:unary}, there must exist graph homomorphisms $\phi_1\colon G \to H$ and $\phi_2\colon H \to G$. In particular, for any $v$ in $V_G$, we must have $\psi_1(v_1) = \phi_1(v)_1$ and $\psi_1(v_2) = \phi_1(v)_2$. Similarly, for any $(u,v)$ in $E_G$, we must have $\psi_1(a_{(u,v)}) = a_{(\phi_1(u),\phi_1(v))}$ and $\psi_1(b_{(u,v)}) = b_{(\phi_1(u),\phi_1(v))}$. Correspondingly, for any $v$ in $V_H$ or $(u,v)$ in $E_H$, we have $\psi_2(v_1) = \phi_2(v)_1$, $\psi_2(v_2) = \phi_2(v)_2$, $\psi_2(a_{(u,v)}) = a_{(\phi_2(u),\phi_2(v))}$, and $\psi_2(b_{(u,v)}) = b_{(\phi_2(u),\phi_2(v))}$. By an argument similar to the preceding one, $\text{id}_{H}$, the identity map on $H$, must be equal to $\phi_1\phi_2$.
		\end{proof}
		
	\end{thm}
	
	\begin{cor} The Retraction Problem for algebras in a language with at least two unary operations is NP-complete.\qed\end{cor}
	
	However, the previous result also immediately demonstrates that
	
	\begin{cor} The Isomorphism Problem for algebras in a language with at least two unary operations is GI-complete.\qed\end{cor}
	
	The complexity of Graph Isomorphism is currently unknown \cite{GareyJohn}, but any algorithm that could prove two graphs were isomorphic in polynomial time would also solve the Isomorphism Problem for these algebras in polynomial time, and vice versa.
	
	With this, we have completely classified the computational complexity of the Homomorphism Factorization Problem and its variants for algebras in a language with at least two unary operations. In the next section, we will finish the proof of Theorem \ref{thm:General}.
	
\section{Rich Languages with an Operation of Arity at Least Two}\label{sec:Binary}

	%Graph Theoretic non-associative result, then associative semigroup. Ternary? Just as bad. Binary and unary? Just as bad. Unary? Unknown, although G-set known to be P (see next section).

	In this section, we consider the Homomorphism Problem in non-associative settings, and then extend this result to classify the computational complexity of all five major variants of the Homomorphism Factorization Problem. This finishes the proof of Theorem \ref{thm:General} and presents an argument that additional requirements must be made of a given variety for Homomorphism Factorization Problem variants to be in polynomial time.
	
	We begin with a demonstration that the Homomorphism Problem is NP-complete for any algebra with a non-associative binary operation. Specifically, we will show that in this case, any algorithmic solution to the find right-factor problem would necessarily give a solution to the Strong Homomorphism Problem for graphs, and vice versa. Strong Graph Homomorphism is known to be NP-complete \cite{Flum}, and therefore this will provide the desired result.
	
	Any undirected graph, $G = (V_G, E_G)$, with at least two vertices is encoded into an algebra $G^{\ast}$ using the following rules -- for every $v$ in $V_G$, there are two elements, $v_1$ and $v_2$ in $G^{\ast}$; and there are four distinguished elements, $a$, $b$, $c$, and $d$. We then assign to $G^{\ast}$ a single, non-associative binary operation, $\cdot$, where for any distinct $u$, $v$ in $V_G$, $\cdot$ is given by Table \ref{tab:Gast}.
	
	\begin{table}[h]
 		\caption{The operation $\cdot$ of $G^{\ast}$.}
		\label{tab:Gast}
		\begin{tabular}{l|llllllll}
			$\cdot$ & $a$ & $b$ & $c$ & $d$ & $u_1$ & $v_1$ & $u_2$ & $v_2$ \\ \hline
 			$a$ & $b$ & $a$ & $a$ & $a$ & $u_1$ & $v_1$ & $u_2$ & $v_2$ \\
 			$b$ & $a$ & $c$ & $a$ & $a$ & $u_1$ & $v_1$ & $u_2$ & $v_2$ \\
    		$c$ & $a$ & $a$ & $d$ & $a$ & $u_1$ & $v_1$ & $u_2$ & $v_2$ \\
    		$d$ & $a$ & $a$ & $a$ & $a$ & $u_1$ & $v_1$ & $u_2$ & $v_2$ \\
    		$u_1$ & $u_1$ & $u_1$ & $u_1$ & $u_1$ & $d$ & $*$ & $c$ & $d$ \\
    		$v_1$ & $v_1$ & $v_1$ & $v_1$ & $v_1$ & $*$ & $d$ & $d$ & $c$ \\
    		$u_2$ & $u_2$ & $u_2$ & $u_2$ & $u_2$ & $c$ & $d$ & $d$ & $b$ \\
    		$v_2$ & $v_2$ & $v_2$ & $v_2$ & $v_2$ & $d$ & $c$ & $b$ & $d$
 		\end{tabular}
	\end{table}
\noindent Note that $*$ is either $u_1v_1 = v_1u_1 = a$ if $(u,v)$ is in $E_G$, or else $u_1v_1 = v_1u_1 = d$. We may intuitively think of $G^{\ast}$ as encoding a degenerate coloring on edges for a new graph produced by connecting the vertices of $G$ to their corresponding counterparts in the complete graph on $V_G$. An example of this construction on the four-cycle, $C_4$, is shown in Fig. \ref{fig:NonAscC4}.
	
	\begin{figure}[h]
					\begin{tikzpicture}[scale=.5]
  						\node (w1) at (0, 0) [circle, draw] {$w_1$};
  						\node (x1) at (3, 2) [circle, draw] {$x_1$};
  						\node (y1) at (8, 2) [circle, draw] {$y_1$};
  						\node (z1) at (5, 0) [circle, draw] {$z_1$};
  						\draw[thick] (w1) -- (x1);
  						\draw[thick] (x1) -- (y1);
  						\draw[thick] (y1) -- (z1);
  						\draw[thick] (z1) -- (w1);
  						\draw[thick] (11,7) -- node[above, text=black] {$a$} ++ (2,0);
						\node (w2) at (0, 6) [circle, draw] {$w_2$};
  						\node (x2) at (3, 8) [circle, draw] {$x_2$};
  						\node (y2) at (8, 8) [circle, draw] {$y_2$};
  						\node (z2) at (5, 6) [circle, draw] {$z_2$};
  						\draw[thick, dotted] (w2) -- (x2);
  						\draw[thick, dotted] (x2) -- (y2);
  						\draw[thick, dotted] (y2) -- (z2);
  						\draw[thick, dotted] (z2) -- (w2);
  						\draw[thick, dotted] (w2) -- (y2);
  						\draw[thick, dotted] (x2) -- (z2);
  						\draw[thick, dotted] (11,5) -- node[above, text=black] {$b$} ++ (2,0);
  						\draw[thick, dotted] (z2) -- (w2);
  						\draw[thick, dotted] (w2) -- (y2);
  						\draw[thick, dashed] (w1) -- (w2);
  						\draw[thick, dashed] (x1) -- (x2);
  						\draw[thick, dashed] (y1) -- (y2);
  						\draw[thick, dashed] (z1) -- (z2);
  						\draw[thick, dotted] (z2) -- (w2);
  						\draw[thick, dotted] (w2) -- (y2);
  						\draw[thick, dashed] (11,3) -- node[above, text=black] {$c$} ++ (2,0);
					\end{tikzpicture}
		\caption{A visual representation of $C_4^{\ast}$.}
		\label{fig:NonAscC4}		
	\end{figure}
	
	This encoding allows us to move between the Strong Graph Homomorphism Problem and the Homomorphism Problem for non-associative algebras, and gives us the following result.
	
	\begin{thm} Let $G$ and $H$ be graphs with at least two vertices. There exists a strong homomorphism $\phi \colon G \to H$ if and only if there exists a homomorphism $\psi \colon G^{\ast} \to H^{\ast}$.
	
		\begin{proof} Suppose first that there exists a strong homomorphism $\phi \colon G \to H$. We construct a function $\psi\colon G^{\ast} \to H^{\ast}$ based on $\phi$ -- specifically, if for any $v$ in $V_G$ we have $\phi(v)$ in $V_H$, then $\psi(v_1) = \phi(v)_1$, $\psi(v_2) = \phi(v)_2$ and for all $x$ in $\{a,b,c,d\}$, $\psi(x) = x$.
		
		We claim $\psi$ is well-defined by the well-definition of $\phi$ and the construction of $H^{\ast}$. Since $v$ maps to a single element of $H$ under $\phi$, that single element, $\phi(v)$, is in turn associated with exactly two distinct elements $\phi(v)_1$ and $\phi(v)_2$ in $H^{\ast}$.
		
		Next, we claim $\psi$ is a homomorphism. Suppose $u$ and $v$ are distinct elements of $V_G$. If $(u,v)$ is in $E_G$, then $u_1v_1 = v_1u_1 = a$, and since $(\phi(u),\phi(v))$ is in $E_H$, $\psi(u_1)\psi(v_1) = \phi(u)_1\phi(v)_1 = a = \psi(a) = \psi(u_1v_1)$. If $(u,v)$ is not in $E_G$, then $u_1v_1 = v_1u_1 = d$, and since $(\phi(u),\phi(v))$ is also not in $E_H$, $\psi(u_1)\psi(v_1) = \phi(u)_1\phi(v)_1 = d =\psi(d) = \psi(u_1v_1)$. Similarly, $\psi(u_2)\psi(v_2) = \phi(u)_2\phi(v)_2 = b = \psi(b) = \psi(u_2v_2)$, $\psi(u_1)\psi(u_2) = \phi(u)_1\phi(u)_2 = c = \psi(c) = \psi(u_1u_2)$, $\psi(u_1)\psi(u_1) = \phi(u)_1\phi(u)_1 = d = \psi(d) = \psi(u_1u_1)$, and $\psi(u_1)\psi(v_2) = \phi(u)_1\phi(v)_2 = d = \psi(d) = \psi(u_1v_2)$. In addition, for all $x$ in $\{a,b,c,d\}$, $\psi(u_1)\psi(x) = \psi(u_1) = \psi(u_1x)$ and $\psi(u_2)\psi(x) = \psi(u_2) = \psi(u_2x)$.
		
		Suppose now that there exists a homomorphism $\psi \colon G^{\ast} \to H^{\ast}$. We begin by arguing that for any $x$ in $\{a,b,c,d\}$, it must be the case that $\psi(x) = x$. Since $d^2 = a$, $\psi(d^2) = \psi(a) = \psi(d)^2$. Since $\psi(d)^2$ is a square, it must be the case that $\psi(a)$ is in $\{a,b,c,d\}$. Suppose by way of contradiction that $\psi(a) = b$. Then $\psi(b) = \psi(a^2) = \psi(a)\psi(a) = b^2 = c$, $\psi(c) = d$, and $\psi(d) = a$, but then $\psi(ad) = \psi(a) = \psi(a)\psi(d) = a$. By similar arguments showing the inconsistency of $\psi(a) = c$ and $\psi(a) = d$, it must be the case that $\psi(a) = a$. It is then immediate that, in fact, $\psi(x) = x$ for any $x$ in $\{a,b,c,d\}$.
		
		Next, suppose $v$ is in $V_G$. We claim that $\psi(v_1) = u_1$ and $\psi(v_2) = u_2$ for some $u$ in $V_H$. We first argue that $\psi(v_i)$ must go to $u_1$ or $u_2$: suppose by way of contradiction that $\psi(v_i)$ were in $\{a,b,c,d\}$. Then $\psi(v_i)^2 \neq \psi(v_i)$ by construction, so in particular $\psi(v_i)(\psi(v_i)^2) = a$, and $\psi(v_i) = a$. Since $a^2 = b$, $\psi(v_i^2) = \psi(v_i)^2 = b$, but $v_i^2 = a$ or $d$ by construction, and by the previous argument these cannot go to $b$ under $\psi$, a contradiction. So it must be the case that $\psi(v_i)$ goes to $u_1$ or $u_2$. Since $G$ has at least two distinct elements, there exists a $w$ in $V_G$ with $v_2w_2 = b$. In particular, $\psi(v_2)\psi(w_2)= \psi(v_2w_2) = \psi(b) = b$, so $\psi(v_2) = u_2$ for some $u$ in $V_H$. Since $v_1v_2 = c$, $\psi(v_1)\psi(v_2) = \psi(v_1v_2) = \psi(c) = c$. So it must be the case that $\psi(v_1)u_2 = c$, so $\psi(v_1) = u_1$. This completes the proof of the claim.
		
		The preceding argument allows us to define $\phi$ based on $\psi$ -- specifically, let $\phi(v)$ be the element in $V_H$ corresponding to $\psi(v_1)$ and $\psi(v_2)$. $\phi$ thus constructed is a well-defined homomorphism.
		\end{proof}
		
	\end{thm}
	
	\begin{cor} The Homomorphism Problem for non-associative algebras is NP-complete.\qed\end{cor}

	Using arguments similar to those in Section \ref{sec:ThreeUnary}, the preceding result can be extended to classify all major variants of the Homomorphism Factorization Problem.
	
	\begin{cor} The Homomorphism Factorization Problem for non-associative algebras is NP-complete. The Exists Right-Factor, Exists Left-Factor, and Retraction Problem for non-associative algebras is NP-complete; however, the Isomorphism Problem for such algebras is GI-complete.\qed\end{cor}
	
	To complete the proof of Theorem \ref{thm:General}, we must now classify the computational complexity of Homomorphism Factorization Problems for semigroups. Specifically, we will show that for associative semigroups, the Exists Right-Factor Problem is NP-complete, and that this result follows from an encoding process similar in spirit to, yet largely distinct from, the previous two demonstrations. This both extends Theorem \ref{thm:General} to a stronger assumption than ``merely'' rich languages, and provides a negative answer to the problem regarding associative semigroups \cite{JosephVanName} that initially motivated this investigation. We then modify this encoding to show that the Retraction Problem and the Exists Left-Factor Problem for semigroups are also NP-complete.
	
	We begin by developing a method for encoding an arbitrary undirected graph, $G = (V_G, E_G)$, as a semigroup, $X_G$. The universe of $X_G$ consists of an element, $v$, for each $v$ in $V_G$; an element, $\chi_{(u,v)}$, for each $u, v$ in $V_G$ such that $(u,v)$ is not an element of $E_G$ (note that we adopt the convention $\chi_{(u,v)} = \chi_{(v,u)}$, unlike in the unary example); distinct elements $b$, $b^2$, and $c$; and a $0$. We assign to $X_G$ the single binary operation, $\cdot$, given by Table \ref{tab:XG}.
	
\begin{figure}[h]
					\begin{tikzpicture}[scale=.5]
  						\node (w1) at (0, 0) [circle, draw] {$w$};
  						\node (x1) at (3, 2) [circle, draw] {$x$};
  						\node (y1) at (8, 2) [circle, draw] {$y$};
  						\node (z1) at (5, 0) [circle, draw] {$z$};
  						\draw[thick] (w1) -- (x1);
  						\draw[thick] (x1) -- (y1);
  						\draw[thick] (y1) -- (z1);
  						\draw[thick] (z1) -- (w1);
  						\node (b) at (4, 6) [circle, draw] {$b$};
  						\draw[scale = 5.0, thick, dotted] (b) to[in=120,out=60,loop] (b);
  						\draw[thick] (b) -- (w1);
  						\draw[thick] (b) -- (x1);
  						\draw[thick] (b) -- (y1);
  						\draw[thick] (b) -- (z1);
  						\draw[thick] (w1) -- (x1);
  						\draw[thick] (x1) -- (y1);
  						\draw[thick] (y1) -- (z1);
  						\draw[thick] (z1) -- (w1);
  						\draw[thick] (b) -- (w1);
  						\draw[thick] (b) -- (x1);
  						\draw[thick] (b) -- (y1);
  						\draw[thick] (b) -- (z1);
  						\draw[thick] (11,6) -- node[above, text=black] {$c$} ++ (2,0);
  						\draw[scale = 5.0, thick, dotted] (b) to[in=120,out=60,loop] (b);
  						\draw[thick, dotted] (11,3) -- node[above, text=black] {$b^2$} ++ (2,0);
					\end{tikzpicture}
		\caption{A visual representation of $X_{C_4}$.}
		\label{fig:XC4}
	\end{figure}
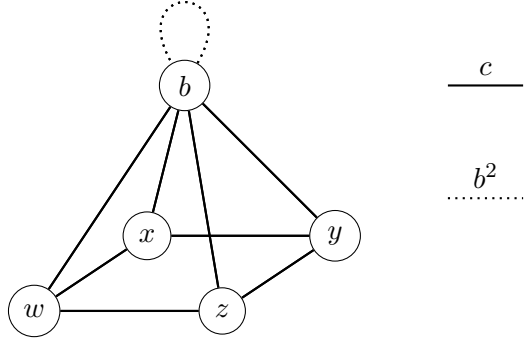

		\begin{table}[h]
 			\caption{The operation $\cdot$ of $X_G$.}
			\label{tab:XG}
			\begin{tabular}{l|lllllll}
				$\cdot$ & $0$ & $b$ & $b^2$ & $c$ & $u$ & $v$ & $\chi$ \\ \hline
 				$0$ & $0$ & $0$ & $0$ & $0$ & $0$ & $0$ & $0$ \\
    			$b$ & $0$ & $b^2$ & $0$ & $0$ & $c$ & $c$ & $0$ \\
    			$b^2$ & $0$ & $0$ & $0$ & $0$ & $0$ & $0$ & $0$ \\
    			$c$ & $0$ & $0$ & $0$ & $0$ & $0$ & $0$ & $0$ \\
    			$u$ & $0$ & $c$ & $0$ & $0$ & $*$ & $*$ & $0$ \\
    			$v$ & $0$ & $c$ & $0$ & $0$ & $*$ & $*$ & $0$ \\
    			$\chi$ & $0$ & $0$ & $0$ & $0$ & $0$ & $0$ & $0$
 			\end{tabular}
		\end{table}
\noindent where for any $u$ and $v$ in $V_G$, $*$ is either $uv = vu = c$ if $(u,v)$ is in $E_G$, or else $uv = vu = \chi_{(u,v)}$; and $\chi$ is a placeholder for any $\chi_{(u,v)}$ in the semigroup. Note that $X_G$ is, in fact, commutative, but this property is not used in the argument. Intuitively, $X_G$ is a description of the graph, $G$, together with a new distinguished vertex, $b$, which is connected to all vertices of $G$. An example of this construction on the four-cycle, $C_4$, is shown in Fig. \ref{fig:XC4}.
	
	We now define $Z$ to be the semigroup with a single binary operation, $\cdot$, given by Table \ref{tab:Zsemigroup}.
		\begin{table}[h]
 			\caption{The operation $\cdot$ of $Z$.}
			\label{tab:Zsemigroup}
			\begin{tabular}{l|lllll}
				$\cdot$ & $0$ & $a$ & $b$ & $b^2$ & $c$ \\ \hline
 				$0$ & $0$ & $0$ & $0$ & $0$ & $0$ \\
    			$a$ & $0$ & $c$ & $c$ & $0$ & $0$ \\
    			$b$ & $0$ & $c$ & $b^2$ & $0$ & $0$ \\
    			$b^2$ & $0$ & $0$ & $0$ & $0$ & $0$ \\
    			$c$ & $0$ & $0$ & $0$ & $0$ & $0$
 			\end{tabular}
		\end{table}
\noindent Note that $Z$ is equivalent to the encoding of the graph consisting of a single loop on a vertex $a$, though it is perhaps more naturally thought of as the encoding of the two element graph which encodes independent sets as homomorphisms \cite{Hell}.
	
	We take as our input two undirected graphs, $G = (V_G, E_G)$ and $H = (V_H, E_H)$. We encode $G$ and $H$ into semigroups, $X_G$ and $Y_H$, using the methods previously discussed. We then construct surjective homomorphisms $f\colon X_G \to Z$ and $h\colon Y_H \to Z$ by taking $f(0) = h(0) = 0$, $f(b) = h(b) = b$, $f(b^2) = h(b^2) = b^2$, and for any $u$ in $V_G$ or $v$ in $V_H$, $f(u) = h(v) = a$, with all other elements going to $c$. This construction is then used in the following theorem.
		
	\begin{thm}There exists a homomorphism $g\colon X_G \to Y_H$ with $f = hg$ if and only if there exists a homomorphism $\phi\colon G \to H$.\label{thm:MainNP} 
		
		\begin{proof} Suppose first that there exists a homomorphism $\phi\colon G \to H$. We construct a function $g\colon X_G \to Y_H$ based on $\phi$ -- specifically, if for any $v$ in $V_G$ we have $\phi(v)$ in $V_H$, we set $g(v) = \phi(v)$, the element of $Y_H$ corresponding to $\phi(v)$. We then set $g(0) = 0$, $g(b) = b$, $g(b^2) = b^2$, and $g(c) = c$. For any $u$, $v$ in $V_G$ with $(u,v)$ not in $E_G$, we send $\chi_{(u,v)}$ to $c$ if $(\phi(u),\phi(v))$ is in $E_H$ or $\chi_{(\phi(u),\phi(v))}$ if it is not.
		
		Such a $g$ is well-defined by construction. We claim that it is also a homomorphism. Suppose that $u$ and $v$ are elements of $X_G$ corresponding to $u$ and $v$ in $V_G$. It is either the case that $uv=vu=c$ or $uv=vu=\chi_{(u,v)}$. Since $g(c) = c$ and $g(\chi_{(u,v)}) = c$ or $\chi_{(\phi(u),\phi(v))}$, we have $g(u)g(v) = \phi(u)\phi(v) = c$ or $\chi_{(\phi(u),\phi(v))} = g(uv)$. Furthermore, we have that $g(u)g(b) = \phi(u)b = c = g(c) = g(ub)$, $g(u)g(c) = \phi(u)c = 0 = g(0) = g(uc)$, and $g(u)g(b^2) = \phi(u)b^2 = 0 = g(0) = g(ub^2)$. Since $g(\chi) = c$ or $\chi$ for any $\chi$ in $X_G$, $g(u)g(\chi) = \phi(u)g(\chi) = 0 = g(0) = g(u\chi)$. We now consider all other possible forms of $x$ and $y$ in $X_G$ -- in the special case where $x=y=b$, we have that $g(b)g(b) = bb = b^2 = g(b^2) = g(bb)$, and in all other cases we have that $g(x)g(y) = 0 = g(0) = g(xy)$. The claim holds.
		
		$g$ also has the property that $f = hg$ -- this is clear from construction for all elements of $X_G$ aside from those coming from $V_G$. Suppose then that $v$ is in $V_G$. Then $g(v) = \phi(v)$, and thus $h(\phi(v)) = a$. Therefore, $hg(v) = a = f(v)$.
		
		Suppose now that there exists a homomorphism $g\colon X_G \to Y_H$ with $f = hg$. We claim that for any $v$ in $V_G$, $g(v)$ is an element of $Y_H$ corresponding to some vertex in $V_H$, $v_{g(v)}$. Since $v$ is in $V_G$, $f(v) = a$, and since $hg = f$, $h(g(v)) = a$. But the only elements which map to $a$ under $h$ are precisely those elements corresponding to vertices in $V_H$; therefore, the claim holds. We then construct a function, $\phi \colon G \to H$, using the rule that $\phi$ sends a given vertex, $v$, in $V_G$ to the vertex in $V_H$ corresponding to $g(v)$. For simplicity, we refer to such a vertex as $g(v)$ -- that is, $\phi(v) = g(v)$, and $g(v)$ is precisely a vertex in $V_H$.
			
			Since $g$ is a homomorphism, $\phi$ is well-defined by construction. It will therefore suffice to show that if $u$ and $v$ are vertices in $V_G$ and $(u,v)$ is an edge in $E_G$, then $(\phi(u),\phi(v))=(g(u),g(v))$ is an edge in $E_H$. Since $(u,v)$ is an edge in $E_G$, $uv = c$ in $X_G$. Consequently, $f(uv) = f(c) = c$, and similarly $h(g(uv)) = h(g(c)) = c$. It must then be either the case that $g(c) = c$ or $g(c) = \chi_{(u',v')}$ for some $u'$ and $v'$ in $V_H$. However, $g$ is a homomorphism, and $g(ub) = g(u)g(b)=g(c)$. Since $f(b) = h(b) = b$, and $b$ is the only element of $X_G$ that maps to $b$ under $f$ and the only element of $Y_H$ that maps to $b$ under $h$, it must be the case that $g(b) = b$. Consequently, $g(ub) = g(u)b = c$, since $g(u)$ corresponds to a vertex in $V_H$ by previous argument. Therefore, $g(c) = c$, and it follows that $g(uv) = g(u)g(v) = c$. Necessarily, then, $(\phi(u),\phi(v))=(g(u),g(v))$ is an edge in $E_H$. Since $\phi$ is well-defined and preserves the edge relation, it is indeed a homomorphism.
		\end{proof}
		
	\end{thm}
		
	\begin{cor} The Exists Right-Factor Problem for semigroups is NP-complete.\qed\end{cor}
	
	As before, we may also use this encoding to classify other Homomorphism Factorization Problems. We need not modify the encoding for the next result:
	
	\begin{thm} Let $G$ and $H$ be undirected graphs. Then $G$ is isomorphic to a subgraph of $H$ if and only if $X_G$ is a retraction of $X_H$.
	
		\begin{proof} Suppose first that $G$ is isomorphic to a subgraph of $H$. Then, in particular, there exist graph homomorphisms $\phi_h \colon G \to H$ and $\phi_g \colon H \to G$ with $\phi_1\phi_2 = \text{id}_H$, the identity map on $H$. We will construct corresponding homomorphisms $h \colon X_G \to X_H$ and $g\colon X_H \to X_G$ with $hg = \text{id}_{X_H}$, the identity map on $X_H$. We begin by setting $g(c)=h(c)=c$, $g(0)=h(0)=0$, $g(b)=h(b)=b$, and $g(b^2)=h(b^2)=b^2$. It then suffices to describe the behavior of $g$ and $h$ on elements corresponding to members of $V_G$ and the $\chi$ elements.
		
		We define $g$ and $h$ based on $\phi_g$ and $\phi_h$, respectively. Specifically, if $v$ in $X_H$ corresponds to an element $v$ in $V_H$, then $g(v) = \phi_g(v)$, the element of $X_G$ corresponding to $\phi_g(v)$ in $V_G$, and if $u$ in $X_G$ corresponds to an element $u$ in $V_G$, then $h(u) = \phi_h(u)$, the element of $X_H$ corresponding to $\phi_h(u)$ in $V_H$. For any $\chi_{(u,v)}$ in $X_H$, if $(\phi_g(u),\phi_g(v))$ is in $E_G$, set $g(\chi_{(u,v)}) = c$, otherwise set $g(\chi_{(u,v)}) = \chi_{(\phi_g(u),\phi_g(v))}$. Similarly, for any $\chi_{(u,v)}$ in $X_G$, if $(\phi_h(u),\phi_h(v))$ is in $E_H$, set $g(\chi_{(u,v)}) = c$, otherwise set $g(\chi_{(u,v)}) = \chi_{(\phi_h(u),\phi_h(v))}$.
		
		It is clear by the proof of Theorem \ref{thm:MainNP} that $g$ and $h$ are homomorphisms. We claim that $\text{id}_{X_H} = hg$. Let $v$ be the element of $X_H$ corresponding to any vertex $v$ in $V_H$. Then $\text{id}_{X_H}(v) = v = \text{id}_H(v) = \phi_h(\phi_g(v)) = h(g(v))$. Similarly, let $\chi_{(u,v)}$ be an element of $X_H$ -- then $(u,v)$ is not an element of $E_H$, and in particular since $G$ is isomorphic to a subgraph of $H$, $(\phi_g(u),\phi_g(v))$ is not an element of $E_G$. Then $\text{id}_{X_H}(\chi_{(u,v)}) = \chi_{(u,v)} = \chi_{(\text{id}_H(u),\text{id}_H(v))} = \chi_{(\phi_h\phi_g(u),\phi_h\phi_g(v))} = h(\chi_{(phi_g(u),\phi_g(v))})=h(g(\chi_{(u,v)}))$. For all other elements $x$ in $X_H$, $g$ and $h$ were already constructed with the property that $g(x) = h(x) = x$, and thus $h(g(x)) = x$. The claim holds.
		
		Suppose now that $X_G$ is a retraction of $X_H$. Then, in particular, there exist homomorphisms $g \colon X_H \to X_G$ and $h \colon X_G \to X_H$ with $hg = \text{id}_{X_H}$, the identity map on $X_H$. We begin by describing the behavior of $g$ and $h$. It clearly must be the case that $g(0) = h(0) = 0$, and this, in turn, means that $g(b) = h(b) = b$, $g(b^2)=h(b^2)=b^2$, and $g(c) = h(c) = c$. If $v$ is an element of $X_H$ corresponding to some $v$ in $V_H$, then $vb = bv = c$, and therefore $g(vb) = g(v)g(b) = g(c) = c = g(v)b$. By the construction of $X_G$, $g(v)$ corresponds to some element $u$ in $V_G$. Suppose $\chi_{(u,v)}$ is in $X_H$, meaning $(u,v)$ is not in $E_H$ and there exists $u$, $v$ in $X_H$ corresponding to some elements $u$, $v$ in $V_H$. Then since $h(c) = c$, it must be the case that $g(\chi_{(u,v)}) = g(u)g(v)= \chi_{(u',v')}$ for some $u'$, $v'$ in $V_G$ with $(u',v')$ not in $E_G$.
		
		We will construct graph homomorphisms $\phi_g\colon H \to G$ and $\phi_h\colon G \to H$ based on $g$ and $h$, respectively. If $v$ is in $V_H$, then there exists an element of $X_H$, $v$, and in particular $g(v)$ is an element of $X_G$ corresponding to some element, $g(v)$, in $V_G$. We set $\phi_g(v) = g(v)$. Similarly, if $v$ is in $V_G$, then there exists an element of $X_G$, $v$, and in particular $h(v)$ is an element of $X_H$ corresponding to some element, $h(v)$, in $V_H$. We set $\phi_h(v) = h(v)$. We claim $\phi_g$ and $\phi_h$ so constructed are homomorphisms. First, suppose $u$ and $v$ are in $V_H$ with $(u,v)$ in $E_H$. Then $g(u)g(v) = c$, meaning $(g(u),g(v))$ is in $E_G$, and since $g(u) = \phi_g(u)$ and $g(v) = \phi_g(v)$, $(\phi_g(u),\phi_g(v))$ is in $E_G$. Next, suppose $u$ and $v$ are in $V_G$ with $(u,v)$ in $E_G$. Then $h(u)h(v) = c$, meaning $(h(u),h(v))$ is in $E_H$, and since $h(u) = \phi_h(u)$ and $h(v) = \phi_h(v)$, $(\phi_h(u),\phi_h(v))$ is in $E_G$. Since $\phi_g$ and $\phi_h$ preserve the edge relation, they are indeed homomorphisms as claimed. Furthermore, by construction it must be the case that $\phi_h\phi_g = \text{id}_H$, the identity map on $H$, since for any vertex $v$ in $V_H$, $\phi_h(\phi_g(v)) = h(g(v)) = v$. Therefore, $G$ is isomorphic to a subgraph of $H$.
		\end{proof}
		
	\end{thm}
	
	\begin{cor} The Retraction Problem for semigroups is NP-complete.\qed\end{cor}
	
	Again as before, this result implies that the Isomorphism Problem for semigroups is GI-complete. However, this result has already been shown in \cite{Booth}. Furthermore, that the Homomorphism Problem for semigroups can be solved in polynomial time, by choosing the homomorphism sending everything to 0, is readily apparent.
	
	This leaves the Exists Left-Factor Problem. Unlike in our previous arguments, the Exists Left-Factor Problem will require a modified construction, and requires the additional assumption that the graphs to be encoded are connected with at least two elements. Therefore, consider a connected, undirected graph, $G = (V_G, E_G)$. Create a new undirected graph $G'$ by taking $V_{G'} = V_G \cup\{w\}$ for some new distinguished vertex $w$, with $w$ not connected to any vertex in $G$. Then construct $X_{G'}$ as per the previous instructions.
	
	We now define $Z'$ to be the semigroup with a single binary operation, $\cdot$, given by Table \ref{tab:Zprimesemigroup}.
		\begin{table}[h]
 			\caption{The operation $\cdot$ of $Z'$.}
			\label{tab:Zprimesemigroup}
			\begin{tabular}{l|llllll}
				$\cdot$ & $0$ & $a$ & $a^2$ & $b$ & $b^2$ & $c$ \\ \hline
 				$0$ & $0$ & $0$ & $0$ & $0$ & $0$ & $0$ \\
    			$a$ & $0$ & $a^2$ & $0$ & $c$ & $0$ & $0$ \\
    			$a^2$ & $0$ & $0$ & $0$ & $0$ & $0$ & $0$ \\
    			$b$ & $0$ & $c$ & $0$ & $b^2$ & $0$ & $0$ \\
    			$b^2$ & $0$ & $0$ & $0$ & $0$ & $0$ & $0$ \\
    			$c$ & $0$ & $0$ & $0$ & $0$ & $0$ & $0$ \\
 			\end{tabular}
		\end{table}
\noindent Note that $Z'$ is equivalent to the encoding of the graph of a single vertex with no edges.
	
	We now take as our input two connected, undirected graphs with at least two elements, $G = (V_G, E_G)$ and $H = (V_H, E_H)$. We encode $G$ and $H$ into semigroups, $X_{G'}$ and $Y_{H'}$, using the methods previously discussed. We then construct homomorphisms $f\colon Z' \to X_{G'}$ and $g\colon Z' \to Y_{H'}$ by taking $f(0) = g(0) = 0$, $f(b) = g(b) = b$, $f(b^2) = h(b^2) = b^2$, $f(a) = w$ for the $w$ in $G'$, $g(a) = w$ for the $w$ in $H'$, $f(a^2) = \chi_{(w,w)}$ in $X_{G'}$, $g(a^2) = \chi_{(w,w)}$ in $Y_{H'}$, and $f(c) = g(c) = c$. This construction is then used in the following theorem.
		
	\begin{thm}There exists a homomorphism $h\colon Y_{H'} \to X_{G'}$ with $f = hg$ if and only if there exists a homomorphism $\phi\colon H \to G$.
		
		\begin{proof} Suppose first that there exists a homomorphism $\phi\colon H \to G$. We construct a function $h\colon Y_{H'} \to X_{G'}$ based on $\phi$ -- specifically, if for any $v$ in $V_H$ we have $\phi(v)$ in $V_G$, we set $h(v) = \phi(v)$, the element of $X_{G'}$ corresponding to $\phi(v)$. We then set $h(0) = 0$, $h(b) = b$, $h(b^2) = b^2$, $h(w) = w$, $h(\chi_{(w,w)}) = \chi_{(w,w)}$ and $h(c) = c$. For any $u$, $v$ in $V_H$ with $(u,v)$ not in $E_H$, we send $\chi_{(u,v)}$ to $c$ if $(\phi(u),\phi(v))$ is in $E_G$ or $\chi_{(\phi(u),\phi(v))}$ if it is not.
		
		Such an $h$ is well-defined by construction. We claim that it is also a homomorphism. Suppose that $u$ and $v$ are elements of $Y_{H'}$ corresponding to $u$ and $v$ in $V_{H'}$. It is either the case that $uv=vu=c$ or $uv=vu=\chi_{(u,v)}$. Since $h(c) = c$ and $h(\chi_{(u,v)}) = c$ or $\chi_{(\phi(u),\phi(v))}$, we have $h(u)h(v) = \phi(u)\phi(v) = c$ or $\chi_{(\phi(u),\phi(v))} = h(uv)$. Furthermore, we have that $h(u)g(b) = \phi(u)b = c = h(c) = h(ub)$, $h(u)h(c) = \phi(u)c = 0 = h(0) = h(uc)$, and $h(u)h(b^2) = \phi(u)b^2 = 0 = h(0) = h(ub^2)$. Since $h(\chi) = c$ or $\chi$ for any $\chi$ in $Y_{H'}$, $h(u)h(\chi) = \phi(u)h(\chi) = 0 = h(0) = h(u\chi)$. We now consider all other possible forms of $x$ and $y$ in $X_G$ -- in the special case where $x=y=b$, we have that $h(b)h(b) = bb = b^2 = h(b^2) = h(bb)$, and in all other cases we have that $h(x)h(y) = 0 = h(0) = h(xy)$. The claim holds. Furthermore, by construction it is clear that $f=hg$ as desired.
		
		Suppose now that there exists a homomorphism $h\colon Y_{H'} \to X_{G'}$ with $f = hg$. By the composition, it must be the case that $h(0) = 0$, $h(b) = b$, $h(b^2) = b^2$, $h(w) = w$, $h(\chi_{(w,w)}) = \chi_{(w,w)}$ and $h(c) = c$. We claim that for any $v$ in $V_H$, $h(v)$ is an element of $X_{G'}$ corresponding to some vertex in $V_G$, $v_{h(v)}$. Since $v$ is in $V_H$, and $H$ is connected with at least two elements, there exists a $u$ in $V_H$ with $(u,v)$ in $E_H$. So $uv = c$, and therefore $h(uv) = h(c) = c = h(u)h(v)$. Thus either $h(v)$ is in $V_H$ or $h(v) = b$ -- however, the second case cannot occur, as $h(vb) = h(c) = c = h(v)h(b) = h(v)b$, and if $h(v) = b$ then $h(v)b = b^2 \neq c$. The claim holds. We then construct a function, $\phi \colon H \to G$, using the rule that $\phi$ sends a given vertex, $v$, in $V_H$ to the vertex in $V_G$ corresponding to $h(v)$. For simplicity, we refer to such a vertex as $h(v)$ -- that is, $\phi(v) = h(v)$, and $h(v)$ is precisely a vertex in $V_G$.
			
			Since $h$ is a homomorphism, $\phi$ is well-defined by construction. It will therefore suffice to show that if $u$ and $v$ are vertices in $V_H$ and $(u,v)$ is an edge in $E_H$, then $(\phi(u),\phi(v))=(h(u),h(v))$ is an edge in $E_G$. Since $(u,v)$ is an edge in $E_H$, $uv = c$ in $Y_{H'}$. Consequently, $h(u)h(v) = h(uv) = h(c) = c$. Necessarily, then, $(\phi(u),\phi(v))=(h(u),h(v))$ is an edge in $E_G$. Since $\phi$ is well-defined and preserves the edge relation, it is indeed a homomorphism.
		\end{proof}
		
	\end{thm}
		
	\begin{cor} The Exists Left-Factor Problem for semigroups is NP-complete.\qed\end{cor}
	
	Finally, we demonstrate a method by which our preceding arguments can be used to show the NP-completeness of the Homomorphism Factorization Problem for algebras with an operation of arbitrarily large arity.
	
	\begin{thm} Let $G = (V_G, E_G)$ be an arbitrary, undirected graph. Then it is possible to encode $G$ as an algebra, $X_G$, with a single ternary operation, such that the Homomorphism Factorization Problem is NP-complete.
	
		\begin{proof} We follow the semigroup encoding to construct $X_G$. The universe of $X_G$ consists of an element, $v$, for each $v$ in $V_G$; an element, $\chi_{(u,v)}$, for each $u, v$ in $V_G$ such that $(u,v)$ is not an element of $E_G$ (note that we adopt the convention $\chi_{(u,v)} = \chi_{(v,u)}$); distinct elements $b$, $b^2$, and $c$; and a $0$. We assign to $X_G$ the single ternary operation, $t(\cdot, \cdot, \cdot)$, given by $t(x,y,z) = x\cdot y$, where $\cdot$ is precisely the associative binary operation given by Table \ref{tab:XG}. Since this construction simultaneously encodes the semigroup from Theorem \ref{thm:MainNP}, it follows that  the Homomorphism Factorization Problem is NP-complete.
		\end{proof}
		
	\end{thm}
	
	\begin{cor} Let $G = (V_G, E_G)$ be an arbitrary, undirected graph. Then it is possible to encode $G$ as an algebra, $X_G$, with a single associative $n$-ary operation, such that the preceding theorem holds, for any $n \geq 3$. \qed \end{cor}
	
	This concludes the proof of Theorem \ref{thm:General}.

\section{Bounded $f$-Cores}\label{sec:Cores}

	We now explore those varieties, still in rich languages, for which certain instances of the Homomorphism Factorization Problem can be solved in polynomial time. To do this, we will investigate the relationship between the Exists Right-Factor Problem and the Homomorphism Problem when $Z$ is fixed. We will then develop a new classification criterion for varieties based on the properties of their finite algebras within this framework.
	
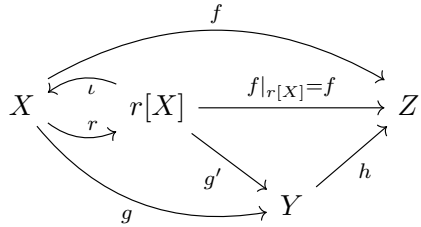
\begin{figure}[h]
			\begin{tikzcd}
				X \arrow[drr, bend right, "g"'] \arrow[rrr, bend left, shift left, shift left, "f"] \arrow[r, leftarrow, bend left, "{\iota}"']  \arrow[r, bend right, "r"] & r[X] \arrow[dr, "{g'}"'] \arrow[rr, "{f\vert_{r[X]} = f}"] && Z\\
				&& Y \arrow[ur, "h"']
			\end{tikzcd}
		\caption{The commutative diagram for Proposition \ref{prp:retract}.}
		\label{fig:retract}
	\end{figure}

	Recall that the general form of Homomorphism Factorization Problems, as shown in Fig. \ref{fig:diagram}, requires as input three finite algebras $X$, $Y$, and $Z$, and a homomorphism, $f\colon X \to Y$, and then asks about the existence of one or both of the homomorphisms $g\colon X \to Y$ and $h\colon Y \to Z$ such that $f =hg$. Suppose that, for a given input, there also exists a retraction $r\colon X \to X$ such that $fr = f$ -- such a retraction will be said to \emph{respect $f$}. We have the following:
	
	\begin{prp}\label{prp:retract} Let $X$, $Y$, and $Z$ be finite algebras, let $f\colon X \to Y$ be a homomorphism, and suppose $r\colon X \to X$ is a retraction that respects $f$. Then $f$ factors through $Y$ if and only if $f \vert_{r[X]}$ factors through $Y$.
		
		\begin{proof} Suppose first that $f$ factors through $Y$. Then there exist homomorphisms $g\colon X \to Y$ and $h\colon Y \to Z$ such that $f =hg$. Since $r$ is a retraction, there exists an inclusion map into $X$, $\iota \colon r[X] \to X$, with $r\iota = r[X]$ and $\iota r = r[X] \subseteq X$. We define a new homomorphism, $g' \colon r[X] \to Y$, by taking $g\iota$ -- since $\iota$ is an inclusion map, and $g$ is a homomorphism, $g'$ is a homomorphism. Furthermore, since $r$ respects $f$, $f\vert_{r[X]} = fr = f = hg = hg\iota = hg'$. Thus, $f \vert_{r[X]}$ factors through $Y$, as desired.
		
		Suppose now that $f \vert_{r[X]}$ factors through $Y$. There there exist homomorphisms $g'\colon r[X] \to Y$ and $h\colon Y \to Z$ such that $f\vert_{r[X]} =hg'$. We define a new homomorphism, $g \colon X \to Y$, by taking $g'r$ -- since $r$ is a retraction, and $g'$ is a homomorphism, $g$ is a homomorphism. Furthermore, since $r$ respects $f$, $f = f\vert_{r[X]}r = hg'r = hg$. Thus, $f$ factors through $Y$, as desired (see Fig. \ref{fig:retract}).
		\end{proof}
	
	\end{prp}
	
	Our classification criterion is related to the existence of retractions that respect $f$. Suppose we have a function $f\colon X \to Z$ for some algebras $X$ and $Z$. We call an algebra $A$ an \emph{$f$-core of an algebra $X$} if $A$ is minimal with respect to the existence of an onto retraction $r\colon X \to A$ that respects $f$. We call $X$ an \emph{$f$-core} if it is its own $f$-core. We classify varieties based on these $f$-cores:
	
	\begin{df} A variety $\mathcal{V}$ has \emph{bounded $f$-cores} if, for any finite algebra $X$ in $\mathcal{V}$ and any surjective homomorphism $f\colon X \to Z$ for which $X$ is an $f$-core, $\vert X \vert \leq s(\vert Z \vert)$ for some function $s$.
	\end{df}
	
	Unlike an arbitrary variety in a rich language, a variety with bounded $f$-cores introduces additional restrictions on Homomorphism Factorization Problems which can reduce their computational complexity. Consider an exists right-factor problem in a variety with bounded $f$-cores where the finite algebra $Z$ is fixed. We have the following:
	
	\begin{thm}\label{thm:fcoreP} The Exists Right-Factor Problem where $Z$ is fixed can be solved in polynomial time for a variety $\mathcal{V}$ if the following conditions hold:
		\begin{enumerate}[I.]
			\item $\mathcal{V}$ has bounded $f$-cores.\\
			\item The $f$-cores of finite algebras in $\mathcal{V}$ can be found in polynomial time.\\
			\item Given a finite algebra in $\mathcal{V}$, a retraction from $X$ to its $f$-core can be found in polynomial time.
		\end{enumerate}
		
		\begin{proof} Suppose $\mathcal{V}$ is a variety for which Conditions I--III hold, and suppose for some fixed algebra $Z$ in $\mathcal{V}$ we have the input of an Exists Right-Factor Problem: algebras $X$ and $Y$, and homomorphisms $f\colon X \to Z$, $h\colon Y \to Z$. We wish to determine whether or not a function $g\colon X \to Y$ with $f = hg$ exists in polynomial time.
		
		Without loss of generality, we may assume $Z = \text{im}(f)$, since it can be checked in polynomial time whether $\text{im}(f) \subseteq \text{im}(h)$, the image of $f$ is known from input, and those elements unique to the image of $h$ may be disregarded. Under this assumption, $f$ is surjective. Furthermore, by Proposition \ref{prp:retract} and Condition III, we may replace $X$ with its $f$-core, $X'$, in polynomial time, and $X'$ can be found in polynomial time by Condition II. By Condition I, there are at most $\vert Y \vert^{\vert X' \vert} = \vert Y \vert^{s(\vert Z \vert)}$ choices for $g'\colon X' \to Y$, and this is a polynomial in $\vert Y \vert$. Our desired homomorphism, then, is given by $g = g'r$, or does not exist, in which case $g'$ will not exist by Proposition \ref{prp:retract}.
 		\end{proof}
	\end{thm}
	
\noindent Some examples of varieties that satisfy all three hypotheses of Theorem \ref{thm:fcoreP} are already known, and will be explored in Section \ref{sec:P}. These examples are especially noteworthy as each of the three conditions required is a nontrivial assumption about the variety in question.

	The first condition of Theorem \ref{thm:fcoreP} alone raises a point of consideration: is there an example of a variety in a rich language without bounded $f$-cores? This is indeed the case; in fact, the variety of semigroups and the variety of semilattices are both varieties for which at least one choice of finite $Z$ can be shown to have unbounded $f$-cores.
	
	\begin{thm}\label{thm:semigroupcore} Let $Z$ be the distinguished semigroup from Section \ref{sec:Binary}. Then for any natural number $n$, there exists a semigroup, $X$, of size at least $n$ that is an $f$-core.
	
		\begin{proof} Recall that we define $Z$ to be the semigroup with a single binary operation, $\cdot$, given by Table \ref{tab:Zsemigroup}.
			
			Let $n$ be a natural number. We encode $K_{n+3}$, the complete graph on $n+3$ vertices, into a semigroup, $X_{K_{n+3}}$, using the method outlined in Section \ref{sec:Binary}. Suppose $f\colon X_{K_{n+3}} \to Z$ is the surjective homomorphism given by taking $f(0) = 0$, $f(b) = b$, $f(b^2) = b^2$, and for any $u$ in $V_{K_{n+3}}$, $f(u) = a$, with all other elements going to $c$. Then by Theorem \ref{thm:MainNP}, any retraction $r\colon X_{K_{n+3}} \to X_{K_{n+3}}$ respecting $f$ has the property of existing if and only if a graph homomorphism $\phi\colon K_{n+3} \to K_{n+3}$ exists. Since any complete graph is a core, such a homomorphism must be the identity map. Consequently, any retraction $r\colon X_{K_{n+3}} \to X_{K_{n+3}}$ respecting $f$ is also the identity map. Therefore, $X_{K_{n+3}}$ is an $f$-core. Furthermore, $\vert X_{K_{n+3}} \vert > n$ by construction.
		\end{proof}	
		
	\end{thm}
	
	\begin{cor} The variety of semigroups does not have bounded $f$-cores. \qed \end{cor}
	
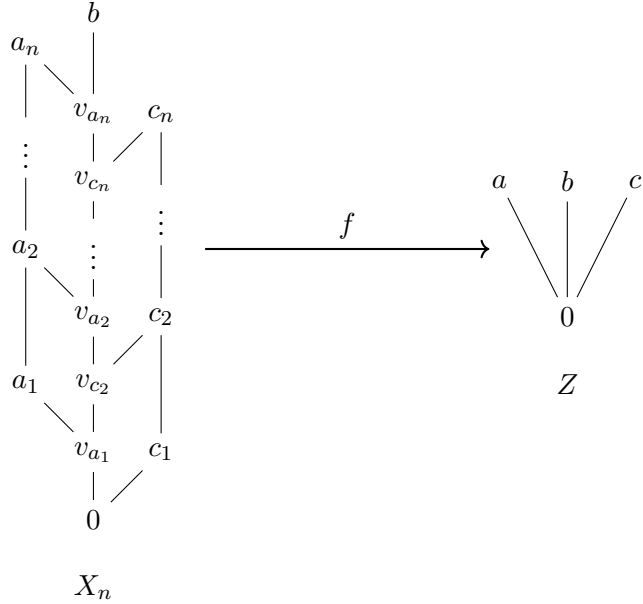
\begin{figure}[h]
			\begin{tikzpicture}[scale=.45]
  				\node (A) at (4, 10) {$a$};
  				\node (B) at (6, 10) {$b$};
  				\node (C) at (8, 10) {$c$};
  				\node (V) at (6, 6) {$0$};
  				\node (S) at (6, 4) {$Z$};
  				\draw (V) -- (A);
  				\draw (V) -- (B);
  				\draw (V) -- (C);
  				\node (start) at (-5, 8) {};
  				\node (end) at (4, 8) {};
				\node (b) at (-8,15) {$b$};
  				\node (a1) at (-10,4) {$a_1$};
  				\node (a2) at (-10,8) {$a_2$};
  				\node (a3) at (-10,11) {$\vdots$};
  				\node (a4) at (-10,14) {$a_n$};
  				\node (va1) at (-8,2) {$v_{a_1}$};
  				\node (vc2) at (-8,4) {$v_{c_2}$};
  				\node (va2) at (-8,6) {$v_{a_2}$};
  				\node (vc3) at (-8,8) {$\vdots$};
  				\node (vc4) at (-8,10) {$v_{c_n}$};
  				\node (va3) at (-8,12) {$v_{a_n}$};
  				\node (c1) at (-6,2) {$c_1$};
  				\node (c2) at (-6,6) {$c_2$};
  				\node (c3) at (-6,9) {$\vdots$};
  				\node (c4) at (-6,12) {$c_n$};
  				\node (zero) at (-8,0) {$0$};
  				\draw (zero) -- (c1) -- (c2) -- (c3) -- (c4);
  				\draw (vc2) -- (c2);
  				\draw (va2) -- (a2);
  				\draw (va3) -- (a4);
  				\draw (vc4) -- (c4);
  				\draw (va1) -- (a1) -- (a2) -- (a3) -- (a4);
  				\draw (zero) -- (va1) -- (vc2) -- (va2) -- (vc3) -- (vc4) -- (va3) -- (b);
  				\draw [thick, ->] (start) -- (end) node[midway, above]{$f$};
  				\node (X) at (-8,-2) {$X_n$};
  			\end{tikzpicture}
		\caption{The construction of an unbounded $f$-core for semilattices.}
		\label{fig:semilattices}
	\end{figure}

	\begin{thm} Consider the semilattice $Z = (\{a,b,c,0\},\wedge)$ given by $a \wedge b = b \wedge c = a \wedge c = 0$. Then for any natural number $n$, there exists a semilattice, $X_n$, of size at least $n$ that is an $f$-core.
	
		\begin{proof} Let $n$ be a natural number -- we create an ``expanded'' version of $Z$, $X_n$, using the following rules: $a$ and $c$ are extended to chains containing $n$ elements, $A_n$ and $C_n$, with labels $a_1$, $a_2$, $\ldots$, $a_n$ and $c_1$, $c_2$, $\ldots$, $c_n$, in ascending order. $0$ is extended to a chain of $2n+1$ elements, $V_n$, with labels $0 = v_{c_1}$, $v_{a_1}$, $v_{c_2}$ , $v_{a_2}$, $\ldots$, $v_{c_n}$, $v_{a_n}$, in ascending order. $b$ is not extended and has the property that $b \geq v$ for any $v$ in $V_n$. We define meet as follows -- for any $i$ from 1 to $n$, $a_i \meet c_i = v_{c_i}$, $a_{i} \meet c_{i+1} = v_{a_i}$, $b \meet a_i = v_{a_i}$, and $b \meet c_i = v_{c_i}$.
			
			$X_n$ so constructed surjects onto $Z$ under the homomorphism $f\colon X_n \to Z$ given by $f[A_n] = a$, $f[V_n] = 0$, $f[C_n] = c$, and $f(b) = b$ (see Fig. \ref{fig:semilattices}). We claim that, in fact, $X_n$ is an $f$-core. Suppose we have a retraction, $e\colon X_n \to X_n$, respecting $f$. By construction, $e[A_n]$ is contained in $A_n$, $e[V_n]$ is contained in $V_n$, $e[C_n]$ is contained in $C_n$, and $e(b) = b$. We will in fact show that, for all possible $x$ in $X_n$, $e(x) = x$.
			
			We first consider the case where $x = a_i$ for some $a_i$ in $A_n$. Now, $v_{a_i} = a_i \meet b$, and thus $e(v_{a_i}) = e(a_i) \meet e(b) = e(a_i) \meet b$. Since $e[A_n]$ is contained in $A_n$, then if $e(x) \neq x$, it must be the case that either $e(a_i) = a_j$ for some $j > i$, or $e(a_i) = a_k$ for some $k < i$. If $e(a_i) = a_j$ for some $j > i$, then $e(v_{a_i}) = a_j \meet b = v_{a_j} > v_{a_i}$. If $j = n$, then since $v_{a_i} < v_{c_n}$, $e(v_{a_i}) < e(v_{c_n})$. So $e(v_{c_n}) > v_{a_n}$, but this is impossible. If $j \neq n$, then similarly since $v_{a_i} < v_{c_{i+1}}$, $e(v_{a_i}) < e(v_{c_{i+1}})$, but by construction $v_{c_{i+1}} < v_{a_j}$, so $e(v_{a_i}) < e(v_{c_{i+1}}) < e(v_{a_j})$. Since $e$ is idempotent, this is impossible. If $e(a_i) = a_k$ for some $k < i$, then $e(v_{a_i}) = a_k \meet b = v_{a_k} < v_{a_i}$. If $k = 1$, then since $v_{a_i} > v_{c_2}$, $e(v_{a_i}) > e(v_{c_2})$. So $e(v_{c_2}) < v_{a_k}$, thus $e(v_{c_2}) = v_{c_1}$. But $v_{a_1} < v_{c_2}$, so $e(v_{a_1}) < v_{c_1}$, and this is impossible. If $k \neq 1$, then similarly since $v_{a_i} > v_{c_{i}}$, $e(v_{a_i}) > e(v_{c_{i}})$, but by construction $v_{c_{i}} > v_{a_k}$, so $e(v_{a_i}) > e(v_{c_{i}}) > e(v_{a_k})$. Since $e$ is idempotent, this is impossible. This, it must be that $e(x) = x$ for this case.
			
			We next consider the case where $x = c_i$ for some $c_i$ in $C_n$. Now, $v_{c_i} = c_i \meet b$, and thus $e(v_{c_i}) = e(c_i) \meet e(b) = e(c_i) \meet b$. Since $e[C_n]$ is contained in $C_n$, then if $e(x) \neq x$, it must be the case that either $e(c_i) = c_j$ for some $j > i$, or $e(c_i) = c_k$ for some $k < i$. If $e(c_i) = c_j$ for some $j > i$, then $e(v_{c_i}) = c_j \meet b = v_{c_j} > v_{c_i}$. If $j = n$, then since $v_{c_i} < v_{a_{n-1}}$, $e(v_{c_i}) < e(v_{a_{n-1}})$. So $e(v_{a_{n-1}}) > v_{c_n}$, but this is impossible. If $j \neq n$, then similarly since $v_{c_i} < v_{a_{i}}$, $e(v_{c_i}) < e(v_{a_i})$, but by construction $v_{a_{i}} < v_{c_j}$, so $e(v_{c_i}) < e(v_{a_{i}}) < e(v_{c_j})$. Since $e$ is idempotent, this is impossible. If $e(c_i) = c_k$ for some $k < i$, then $e(v_{c_i}) = c_k \meet b = v_{c_k} < v_{c_i}$. If $k = 1$, then since $v_{c_i} > v_{a_1}$, $e(v_{c_i}) > e(v_{a_1})$. So $e(v_{a_1}) < v_{c_k} = v_{c_1}$, but this is impossible. If $k \neq 1$, then similarly since $v_{c_i} > v_{a_{i-1}}$, $e(v_{c_i}) > e(v_{a_{i-1}})$, but by construction $v_{a_{i-1}} > v_{c_k}$, so $e(v_{c_i}) > e(v_{a_{i-1}}) > e(v_{c_k})$. Since $e$ is idempotent, this is impossible. This, it must be that $e(x) = x$ for this case.
			
			Lastly, suppose $x = v$ for some $v$ in $V_n$. Without loss of generality, suppose $x = v_{a_i}$ for some $a_i$ in $A_n$. Then $x = a_i \meet b$, and thus $e(x) = e(a_i) \meet b$. By the preceding argument, it must be the case that $e(a_i) = a_i$; therefore, $e(x) = x$.
			
			Therefore, $X_n$ is an $f$-core, and by construction, $\vert X_n \vert > n$.
		\end{proof}	
		
	\end{thm}
	
	\begin{cor} The variety of semilattices does not have bounded $f$-cores. \qed \end{cor}
	
	The second example is especially significant as the semilattice $Z$ is a projective object in the variety of semilattices \cite{Horn}. Thus, we have also demonstrated that projectivity, a condition similar to retraction, is insufficient to guarantee the existence of bounded $f$-cores. It remains an open question precisely what conditions a variety must satisfy to always have bounded $f$-cores.

	Condition II from Theorem \ref{thm:fcoreP}, requiring that $f$-cores can be found in polynomial time, does not hold in general. The previous constructions for semigroups will suffice, as it has been shown in \cite{HellCore} that determining the core of a graph is NP-complete.
	
	\begin{thm} Any algorithm that can find the $f$-core of an arbitrary semigroup, $X$, is capable of finding the minimal retraction of an arbitrary graph, $G = (V_G, E_G)$.
		
		\begin{proof} We consider the special case used in Theorem \ref{thm:semigroupcore} -- let $X$ be the encoding of an arbitrary graph $G$, and $Z$ be the distinguished semigroup. Let $f \colon X \to Z$ be the surjective homomorphism given by taking $f(0) = 0$, $f(b) = b$, $f(b^2) = b^2$, and for any $u$ in $V_{G}$, $f(u) = a$, with all other elements going to $c$. Then, if $X'$ is the $f$-core of $X$, there exists a retraction $r \colon X \to X$ with $r[X] = X'$ and $fr = f$. Furthermore, $r'[X'] = X'$ for all other retracts $r' \colon X \to X$ by minimality. However, by Theorem \ref{thm:MainNP}, $r$ exists if and only if a graph homomorphism $\phi_r \colon G \to G$ exists. Furthermore, $X' = r[X]$ corresponds to a subgraph of $X$, and by the minimality of $X'$, this subgraph is also minimal. Such a subgraph is the minimal retraction of $G$.
		\end{proof}
		
	\end{thm}
	
	\begin{cor} The problem of finding the $f$-core of an arbitrary algebra is NP-complete. \qed \end{cor}

\noindent It should be noted that semigroups are already known to not have bounded $f$-cores, and that for every variety currently known to have bounded $f$-cores, the $f$-cores of arbitrary finite algebras can be found in polynomial time. Consequently, it may very well be the case that Condition II is an immediate consequence of Condition I.

	Condition III similarly does not hold in general. Once again, the previous constructions for semigroups will provide a counterexample.
	
	\begin{thm} Any algorithm that can produce a function mapping an arbitrary semigroup, $X$, to its known $f$-core is capable of mapping an arbitrary graph $G = (V_G, E_G)$ to its known core.
		
		\begin{proof} Suppose $G'$ is known to be the core of an arbitrary graph $G$. If $X_{G'}$ and $X_G$ are the encoding, into semigroups, of $G'$ and $G$, respectively, $Z$ is the distinguished semigroup, and $f \colon X \to Z$ is the surjective homomorphism given by taking $f(0) = 0$, $f(b) = b$, $f(b^2) = b^2$, and for any $u$ in $V_{G}$, $f(u) = a$, with all other elements going to $c$, we claim that then $X_{G'}$ is the $f$-core of $X_G$. Since $G'$ is a core, $X_{G'}$ is an $f$-core by minimality -- furthermore, $X_{G'}$ is contained in $X_G$ under the inclusion map induced by the inclusion of $G'$ in $G$ under Theorem \ref{thm:MainNP}. The claim holds. Therefore, also by Theorem \ref{thm:MainNP}, any map from $X_G$ to $X_{G'}$ induces a map from $G$ to $G'$.
		\end{proof}
		
	\end{thm}
	
	\begin{cor} The problem of mapping an arbitrary algebra to its $f$-core is NP-complete. \qed \end{cor}
	
\noindent As before, it should be noted that semigroups are already known to not have bounded $f$-cores, and that for every variety currently known to have bounded $f$-cores, the $f$-cores of arbitrary finite algebras can be mapped to in polynomial time. Consequently, it may very well be the case that Condition III is an immediate consequence of Condition II or Condition I.

	The previous results suggest an underlying relationship between the three conditions in Theorem \ref{thm:fcoreP}. This, in turn, suggests that the notion of bounded $f$-cores is tied to polynomial instances of the Exists Right-Factor Problem. We speculate that the following may be the complete classification of the Exists Right-Factor Problem.
	
	\begin{con}\label{conjecture} The Exists Right-Factor Problem with fixed algebra $Z$ is in polynomial time for a given variety if and only if the variety has bounded $f$-cores.
	\end{con}

\noindent We suspect that Conjecture \ref{conjecture} is the case as there are both no known examples of varieties without bounded $f$-cores for which this version of the Exists Right-Factor Problem is known to always be in polynomial time. In addition, it appears that conditions I through III cannot occur independently of one another. Finally, every known variety for which these instances of the Exists Right-Factor Problem are always in polynomial time.

\section{Polynomial Time Results}\label{sec:P}

	In this section, we explore those special cases for which the Exists Right-Factor Problem with fixed algebra $Z$ can always be solved in polynomial time. In particular, these varieties all satisfy all three hypotheses of Theorem \ref{thm:fcoreP}, and thus serve as evidence in support of Conjecture \ref{conjecture}. We believe that these varieties may also suggest techniques which can be used to produce polynomial time algorithms for additional examples.
	
	We begin by demonstrating that, for the first three tame congruence types, the Exists Right-Factor Problem with fixed algebra $Z$ is in polynomial time. Currently, the computational complexity of these problems for semilattices and lattices is unknown; however, we suspect that the property of having bounded $f$-cores is not related to congruence type.
	
	\begin{thm} Let $\mathcal V$ be the variety of $G$-sets for some finite group $G$. Then Conditions I-III of Theorem \ref{thm:fcoreP} hold.
	
		\begin{proof} Fix a finite group $G$, and let $X$ and $Z$ be finite $G$-sets with some surjective homomorphism $f\colon X \to Z$. Since $f$ is a $G$-set homomorphism, we have, in particular, that $gf(x) = f(gx)$ for any $g \in G$, $x \in X$. Suppose, then, that $r \colon X \to X$ is a retraction that respects $f$ -- that is, $fr = f$. Since $r$ is a retraction, it must preserve the orbits of $X$ under the action from $G$. In particular, we may take $r$ to send each element of an orbit to a distinguished representative of that orbit -- since the orbits partition $X$, such a definition is well-defined. Furthermore, $r$ can be produced in polynomial time, since $G$ is fixed and the orbits can therefore be calculated in time polynomial to $\vert X \vert$. $r$ is also minimal by construction.
		
		We claim that the resulting $f$-core, $X'$, is bounded. Since $f$ is a $G$-set homomorphism, $f$ sends orbits to corresponding orbits in $Z$. By construction, $\vert X' \vert$ is equal to the number of orbits of $Z$ by the surjectivity of $X$. Furthermore, by the class equation of the group action of $G$ on $Z$, $\vert Z \vert = \sum_{i \in \mathcal{I}} \vert Gz_i \vert$ where, for each $i \in \mathcal{I}$, $Gz_i$ is a distinct orbit of $Z$. Therefore, $\vert X' \vert \leq s(\vert Z \vert)$ for some function on $Z$, and therefore the claim holds.
		
		We have therefore constructed and mapped to a bounded $f$-core in polynomial time, satisfying all three conditions of Theorem \ref{thm:fcoreP}.
		\end{proof}	
		
	\end{thm}
	
	\begin{cor} Let $\mathcal V$ be the variety of $G$-sets for some finite group $G$. Then for any finite $Z$ in $\mathcal{V}$, the Exists Right-Factor Problem with fixed algebra $Z$ is in polynomial time. \qed \end{cor}
		
\noindent Interestingly, while the computational complexity for the general case of Homomorphism Factorization Problems in languages with a single unary operations remains unknown, the above corollary provides us with at least one classification of a variety that is not in a rich language.

	Similarly, vector spaces and Boolean algebras can also be demonstrated to be definitively polynomial time examples of this problem using the following argument.
	
	\begin{thm} Let $\mathcal V$ be the variety of vector spaces over a field $F$. Then Conditions I-III of Theorem \ref{thm:fcoreP} hold.
	
		\begin{proof} In the case where $F$ is infinite, this will be trivially true as all finite vector spaces will be isomorphic to the $0$-vector space. Fix $F$ to be a finite field, and let $X$ and $Z$ be finite vector spaces over $F$ with $f\colon X \to Z$ a surjective homomorphism. By the rank-nullity theorem, $\dim (\ker f) + \dim (\text{im }f) = \dim X$ \cite{Roman}, and since $f$ is surjective, $\dim (\text{im }f) = \dim Z$. Choose $r \colon X \to X$ to be the retraction sending $X$ to $X' = X/\ker f$ chosen by selecting a distinguished element of each conjugacy class. This selection can be performed in polynomial time using the data for $f$, and $X'$ is an $f$-core by minimality. Since $r$ is a homomorphism, we have that $\dim (\ker r) + \dim (\text{im }r) = \dim X$. Since $fr =f$, $\dim (\text{im }r) \geq \dim (\text{im }f)$, and therefore $\dim (\ker r) \leq \dim (\ker f)$. Therefore, $\dim (\text{im }r) \leq \dim X - \dim (\ker f) = \dim (\text{im }f)$. Since $X$ and $Z$ are finite, it follows that $\vert X' \vert \leq \vert Z \vert$.
		
		We have therefore constructed and mapped to a bounded $f$-core in polynomial time, satisfying all three conditions of Theorem \ref{thm:fcoreP}.
		\end{proof}	
		
	\end{thm}
	
	\begin{cor} Let $\mathcal V$ be the variety of vector spaces over a field $F$. Then for any finite $Z$ in $\mathcal{V}$, the Exists Right-Factor Problem with fixed algebra $Z$ is in polynomial time. \qed \end{cor}
	
	\begin{thm} Let $\mathcal V$ be the variety of Boolean algebras. Then Conditions I-III of Theorem \ref{thm:fcoreP} hold.
	
		\begin{proof} Let $X$ and $Z$ be finite Boolean algebras with $f\colon X \to Z$ a surjective homomorphism. Since $f$ is surjective, and $X$ and $Z$ are finite, it will suffice to consider the atoms of the Boolean algebras \cite{Givant}. We construct $r \colon X \to X$ with $fr=f$ by sending each atom to a distinguished representative of its image under $f$ -- this construction can clearly be performed in polynomial time. Then, in fact, we have $r[X] = X' \cong Z$ by the classification of finite Boolean algebras, and hence $X'$ is necessarily an $f$-core which is bounded.
		
		We have therefore constructed and mapped to a bounded $f$-core in polynomial time, satisfying all three conditions of Theorem \ref{thm:fcoreP}.
		\end{proof}	
		
	\end{thm}
	
	\begin{cor} Let $\mathcal V$ be the variety of Boolean algebras. Then for any finite $Z$ in $\mathcal{V}$, the Exists Right-Factor Problem with fixed algebra $Z$ is in polynomial time. \qed \end{cor}

\noindent These examples also provide instances of varieties in rich languages for which the Exists Right-Factor Problem with fixed algebra $Z$ is in polynomial time regardless of the chosen $Z$. This suggests the existence of a property or identity, currently unknown, which distinguishes these varieties from those for which known NP-complete instances exist.
	
	The potential nature of such a property was explored using a simpler construction, still in a rich language. This also provides us with an example of a variety for which the Exists Right-Factor Problem with fixed algebra $Z$ is in polynomial time but is not a tame congruence type.
	
	\begin{thm} Let $\mathcal V$ be the variety of abelian groups. Then Conditions I-III of Theorem \ref{thm:fcoreP} hold.
	
		\begin{proof} Let $X$ and $Z$ be finite abelian groups with $f\colon X \to Z$ a surjective homomorphism. Since $X$ and $Z$ are finite, by the Primary Decomposition Theorem for finite abelian groups both $X$ and $Z$ are isomorphic to direct products of their Sylow subgroups, which are cyclic \cite{Dummit}. It is therefore sufficient to consider the case where both $X$ and $Z$ are cyclic, as our retraction will consist of the componentwise retraction on the cyclic portions. Since $f$ is a surjection, $\vert Z \vert \vert \vert X \vert$, and indeed there is an isomorphic copy of $Z$ in $X$. Consequently, we let $r \colon X \to X$ be the retraction to this copy of $Z$ -- $r$ can be constructed in polynomial time, and $r[X] = X'$ is the $f$-core of $X$ by construction with, indeed, $\vert X' \vert = \vert Z \vert$.
		
		We have therefore constructed and mapped to a bounded $f$-core in polynomial time, satisfying all three conditions of Theorem \ref{thm:fcoreP}.
		\end{proof}	
		
	\end{thm}
	
	\begin{cor} Let $\mathcal V$ be the variety of abelian groups. Then for any finite $Z$ in $\mathcal{V}$, the Exists Right-Factor Problem with fixed algebra $Z$ is in polynomial time. \qed \end{cor}
	
\section{Conclusions}\label{sec:conclusions}

	A characterization of varieties with bounded $f$-cores appears to be the most productive route to completing the classification of Homomorphism Factorization Problems at this time. The relationship between $f$-cores and graph cores already apparent would seem to suggest that many graph-theoretic results could have applications to algebraic computational questions which have not previously been addressed -- the result of Hell and Ne\v set\v ril linking graph core homomorphisms to graph homomorphisms \cite{HellCore} itself suggests that $f$-cores may play a fundamental role in the nature of algebraic homomorphisms.
	
	Further evidence for or against Conjecture \ref{conjecture} most likely regards the variety of semilattices -- it is currently unknown whether or not the Exists Right-Factor Problem for semilattices is NP-complete, though we suspect this is indeed the case. This, in addition to a proof regarding any possible relationship between Conditions I-III of Theorem \ref{thm:fcoreP}, seems to be the most approachable method of tackling the conjecture at this time. It is not known what a general proof of the conjecture would necessarily include, but may depend on the characterization of bounded $f$-cores.


\begin{bibdiv}
\begin{biblist}
\bib{bimmvw}{article}{
   author={Berman, Joel},
   author={Idziak, Pawe\l },
   author={Markovi\'c, Petar},
   author={McKenzie, Ralph},
   author={Valeriote, Matthew},
   author={Willard, Ross},
   title={Varieties with few subalgebras of powers},
   journal={Trans. Amer. Math. Soc.},
   volume={362},
   date={2010},
   number={3},
   pages={1445--1473},
   issn={0002-9947},
   review={\MR{2563736}},
}
\bib{Booth}{article}{
	author = {Booth, Kellogg S.},
	title = {Isomorphism testing for graphs, semigroups, and finite automata are polynomially equivalent problems},
	journal = {SIAM Journal on Computing},
	volume = {7},
	number = {3},
	pages = {273-279},
	year = {1978},
	doi = {10.1137/0207023},
}
\bib{Dummit}{book}{
   author={Dummit, David S.},
   author={Foote, Richard M.},
   title={Abstract algebra},
   edition={3},
   publisher={John Wiley \& Sons, Inc., Hoboken, NJ},
   date={2004},
   pages={xii+932},
   isbn={0-471-43334-9},
   review={\MR{2286236}},
}
\bib{Flum}{book}{
   author={Flum, J\"{o}rg},
   author={Grohe, Martin},
   title={Parameterized Complexity Theory},
   publisher={Springer-Verlag Berlin Heidelberg, Germany},
   date={2006},
   pages={xiii+495},
   isbn={978-3-540-29953-0},
}
\bib{GareyJohn}{book}{
   author={Garey, Michael R.},
   author={Johnson, David S.},
   title={Computers and intractability},
   note={A guide to the theory of NP-completeness;
   A Series of Books in the Mathematical Sciences},
   publisher={W. H. Freeman and Co., San Francisco, Calif.},
   date={1979},
   pages={x+338},
   isbn={0-7167-1045-5},
   review={\MR{519066}},
}
\bib{Givant}{book}{
   author={Givant, Steven},
   author={Halmos, Paul},
   title={Introduction to Boolean algebras},
   series={Undergraduate Texts in Mathematics},
   publisher={Springer, New York},
   date={2009},
   pages={xiv+574},
   isbn={978-0-387-40293-2},
   review={\MR{2466574}},
}
\bib{Hell}{book}{
   author={Hell, Pavol},
   author={Ne\v set\v ril, Jaroslav},
   title={Graphs and homomorphisms},
   series={Oxford Lecture Series in Mathematics and its Applications},
   volume={28},
   publisher={Oxford University Press, Oxford},
   date={2004},
   pages={xii+244},
   isbn={0-19-852817-5},
   review={\MR{2089014}},
}

\bib{HellCore}{article}{
   author={Hell, Pavol},
   author={Ne\v set\v ril, Jaroslav},
   title={The core of a graph},
   note={Algebraic graph theory (Leibnitz, 1989)},
   journal={Discrete Math.},
   volume={109},
   date={1992},
   number={1-3},
   pages={117--126},
   issn={0012-365X},
   review={\MR{1192374}},
}
\bib{Horn}{article}{
   author={Horn, Alfred},
   author={Kimura, Naoki},
   title={The category of semilattices},
   journal={Algebra Universalis},
   volume={1},
   date={1971},
   number={1},
   pages={26--38},
   issn={0002-5240},
   review={\MR{0318019}},
}
\bib{Roman}{book}{
   author={Roman, Steven},
   title={Advanced linear algebra},
   series={Graduate Texts in Mathematics},
   edition={3},
   publisher={Springer, New York, NY},
   date={2008},
   pages={xviii+525},
   isbn={978-0-387-72828-5},
}
\bib{JosephVanName}{webpage}{
	author = {Van Name, J.},
	url = {https://mathoverflow.net/questions/261966/},
	title = {Is there a good computer program for searching for endomorphisms between finite algebras which make diagrams commute? Is this problem NP-complete?},
	date = {Feb 11 '17 at 18:25},
	accessdate = {Feb 14 '18},
}
\end{biblist}
\end{bibdiv}
\end{document}